\documentclass{elsart3-1}

\usepackage[active]{srcltx}
\usepackage{amssymb}
\usepackage{amsmath}
\usepackage[frenchb]{babel}
\usepackage[applemac]{inputenc}
\usepackage{enumerate}
\newtheorem{e-proposition}[theorem]{Proposition}

\newtheorem{e-definition}[theorem]{Definition\rm}

\newtheorem{theoreme}{Th\'eor\`eme}

\def\og{\leavevmode\raise.3ex\hbox{$\scriptscriptstyle\langle\!\langle$~}}
\def\fg{\leavevmode\raise.3ex\hbox{~$\!\scriptscriptstyle\,\rangle\!\rangle$}}

\newcounter{monPar}
\newcommand{\parag}{\refstepcounter{monPar}%
\par\smallskip\noindent\textbf{\themonPar.} }

\def \RM{\mathbb{R}}
\def \ZM{\mathbb{Z}}%

\def \Vol {{\rm Vol\,}}
\def \d{\partial}
 
\def\a{\alpha}

\def\p{\varphi}

\def \s{\sigma}

\def \to{\longrightarrow} 
\def \w{\wedge}

\def \< {{\langle }}
\def \> {{\rangle }}

\begin{document}
\centerline{Syst{\`e}mes dynamiques}
\title{Sur le th\'eor\`eme KAM.}
\author{Mauricio Garay}
\ead{garay91@gmail.com}
\address{ Institut für Mathematik\\
FB 08 - Physik, Mathematik und Informatik\\
Johannes Gutenberg-Universität Mainz\\
Staudinger Weg 9\\
55128 Mainz.}

\medskip
\selectlanguage{francais}

\noindent{\bf Résumé.}
\vskip 0.5\baselineskip
\noindent
Le théorème KAM garantit la stabilité, sous l'effet de perturbations, de certains tores invariants dans les systèmes hamiltoniens.  Cependant l'hypothèse de non-dégénérescence de l'application des fréquences, nécessaire à l'application de ce théorème, n'est pas directement vérifiée en mécanique céleste. Cette difficulté conduit souvent à des calculs fastidieux et difficiles. Dans cette note, j'annonce  des théorèmes de tores invariants qui évitent cette hypothèse et, par conséquent, de tels calculs. Ces résultats permettent de répondre affirmativement à la conjecture de Herman sur les tores invariants des symplectomorphismes analytiques.
\vskip 0.5\baselineskip

\noindent{\bf Abstract.}
\vskip 0.5\baselineskip
\noindent
{\bf On the KAM theorem. }
The KAM theorem ensures the persistence, under perturbations, of some special invariant tori in hamiltonian dynamical systems. However, the non-degeneracy hypothesis of the frequency map, which is necessary to apply the theorem, is not fullfiled in celestial mechanics. This problem gives rise to complicated and tedious computations. I announce here invariant tori theorems which avoid this assumption and therefore such computations.
These results give a positive answer to Herman's invariant tori conjecture for analytic symplectomorphism.

\bigskip

\parag Le théorème KAM fournit un critère de persistence des tores invariants d’un système dynamique hamiltonien.  Dans la version  initiale de Kolmogorov, deux conditions étaient requises pour pouvoir appliquer le théorème~\cite{Kolmogorov_KAM}. La  première~–~appelée condition  diophantienne~–~porte sur  la dynamique le long du tore invariant initial. La seconde~–~dite de non-dégénérescence isochronique~–~demande qu'une certaine application des périodes soit  un isomorphisme. Cette deuxième condition a été affaiblie par Rüssmann, et on pensait que son résultat 
était optimal~\cite{Russmann_KAM,Russmann_degeneracy}. Nous verrons que ce n’est pas exactement le cas.  

\parag Commençons par  rappeler le contexte du théorème KAM, dans sa version réelle analytique. Soit donc $U$ un ouvert connexe de $T^*\RM^n=\RM^n \times \RM^n$ 
muni des coordonnées $q_1,...,q_n, p_1,...,p_n $, ainsi que de sa structure symplectique habituelle $\omega=\sum_{i=1}^ndq_i \w dp_i $. 
Une fonction analytique réelle
$$H: U \to \RM$$
déﬁnit un {\em système dynamique Hamiltonien} donné par les équations de Hamilton~: 
$$ \dot q_i= \d_{p_i}H,\  \dot p_i=-\d_{q_i} H .$$
 
Ce système dynamique est appelé {\em Liouville-intégrable}, ou tout simplement {\em intégrable}, s’il existe des fonctions analytiques
$f_1,...,f_n$ dont les crochets de Poisson s’annulent deux à deux,  et si les ﬁbres de {\em l’application moment}~: 
$$f=(f_1,...,f_n): U \to \RM^n $$
sont de dimension $n$. Si l’application moment est propre et lisse, les ﬂots des $f_i$ induisent une structure localement affine sur 
les ﬁbres de cette application $f$, ce sont donc des tores. Dans ce cas, le théorème d’existence des coordonnées
actions-angles affirme que l'on peut linéariser localement l'application moment et la dynamique par un changement de variables symplectique~: quitte à se restreindre à un voisinage suffisamment petit d'un tore, on peut trouver des coordonnées $\p_i \in (\RM/2\pi \ZM),I_i \in \RM$, $i=1,\dots,n$, dans lesquelles
\begin{enumerate}
\item $U$ est un voisinage de la section nulle du ﬁbré cotangent  $T^*(S^1)^n=\{ (\p,I) \}$ définit comme préimage par $f$ d'un voisinage de l'origine  ;
\item $f_i(\p,I)=I_i$ et $\omega={\displaystyle \sum_{i=1}^n} d\p_i  \w dI_i $ ;
\end{enumerate}
où l'on a posé $\p=(\p_1,\p_2,\dots,\p_n)$ et $I=(I_1,I_2,\dots,I_n)$. La fonction $H$ commute avec les $f_i$, elle ne dépend donc que de la valeur de $I$ et pas de celle de $\p$. Dans ces nouvelles coordonnées,  les équations de Hamilton s'écrivent sous la forme~:
$$\dot \p_i= \d_{I_i}H,\ \dot I_i=0. $$
 Le ﬂot de $H$ est alors linéaire sur les fibres de l'application $f$~\cite{Arnold_Liouville,Mineur}.
Un tel système dynamique sur le tore est appelé {\em quasi-périodique}. Dans les coordonnées $\p,I$, il est
déterminé par le vecteur vitesse des trajectoires, appelé le {\em vecteur des fréquences.} Pour un système intégrable propre, on définit ainsi {\em l'application des fréquences}
$$F:U \to \RM^n ,$$
constante sur les fibres de $f$, et qui associe à chaque tore le vecteur de fréquences correspondant. Dans les coordonnées actions-angles, l'application $F$ est simplement le gradient de l'application $H$, comme le montrent les équations de Hamilton.
 
\parag Nous dirons qu’un vecteur 
$\omega=(\omega_1,...,\omega_n) \in \RM^n$
est {\em diophantien} s'il existe des constantes $(C,\tau)$ telles que : 
$$ \forall i \in \ZM^n \setminus \{ 0 \},\  |(\omega,i)| ≥\frac{C}{\| i\|^\tau}.$$

Un {\em tore KAM} d'un système hamiltonien est un tore invariant sur lequel la dynamique est quasi-périodique de fréquence diophantienne.
Le théorème KAM affirme que pour des perturbations suffisamment petites d’un système intégrable propre et lisse,
il existe un ensemble de mesure non vide, réunion de tores KAM, {\em pourvu que l'application des périodes soit de rang maximal}. Cette dernière condition est appelée la {\em non-dégénérescence isochoronique}. Dans la pratique, les symétries imposées par la physique entraînent que les systèmes hamiltoniens peuvent être fortement dégénérés, comme c'est le cas en mécanique
céleste. La transformation du système en un système non-dégénéré nécessite des calculs et des artifices compliqués, y compris pour le problème des trois corps restreint~\cite{Arnold_trois_corps}. Pour ces raisons, plusieurs tentatives ont été faites pour affaiblir l'hypothèse de non-dégénérescence.
 
Le résultat le plus aboutit, dans cette direction, est celui de Rüssmann,  qui a montré un théorème de tores invariants sous l'hypothèse que l'image de l’application des  fréquences n'est pas contenue dans un hyperplan~\cite{Russmann_KAM}.  

\parag  Reprenant une construction de Katok, Sevryuk a montré que la condition de Rüssmann ne pouvait être évitée~\cite{Katok_KAM,Sevryuk_KAM}.   Il existe cependant de nombreux tores invariants même dans les cas dégénérés, comme le montre le résultat suivant~: 

\begin{theoreme} Le voisinage d'un tore KAM d'un système dynamique hamiltonien analytique possède un ensemble de mesure positive de tores KAM.
De plus, la densité de cet ensemble est égale à un au voisinage de chaque tore KAM.
\end{theoreme}
Rappelons que la {\em densité} d'un ensemble mesurable $\Omega \subset \RM^n$ en un point $x \in \RM^n$, lorsqu'elle existe, est définie par la formule
$$\lim_{r \to 0} \frac{\Vol(B(x,r) \cap \Omega)}{\Vol(B(x,r))} $$
où $\Vol(-)$ désigne la mesure de Lebesgue et $B(x,r)$ la boule de centre $x$ et de rayon $r$.

Herman avait demandé si un système hamiltonien peut posséder des tores KAM isolés, ce théorème montre que c'est impossible.

On peut aller plus loin~: le même résultat reste valable lorsque l'on introduit des paramètres. Ainsi, contrairement à ce que l'on pourrait penser, toute perturbation d'un système intégrable possède de nombreux tores KAM, même lorsque l'hypothèse de non-dégénérescence de Rüssmann n'est pas vérifiée.  

\parag Le dernier résultat que je viens d'évoquer peut sembler en contradiction avec les exemples de Katok et Sevryuk. Considérons, avec ces auteurs, le cas deux oscillateurs harmoniques couplés sans interactions $H=I_1+\sqrt{2}I_2 $. Choisissons une approximation rationnelle $p/q$ de $\sqrt{2}$ et posons
$$H'=H+(\frac{p}{q}-\sqrt{2})I_2+\cos(p\p_1-q\p_2). $$
Le nouveau système hamiltonien $H'$ possède alors deux intégrales premières $\cos(p\p_1-q\p_2),\ qI_1+p I_2 $, dont les niveaux définissent des cylindres. On vérifie que presque tous les mouvements décrivent des hélices tracées sur ces cylindres. Celles-ci sont non-bornées, elles ne peuvent donc  être confinées sur des tores. En prenant des approximations rationnelles de plus en plus précises de $\sqrt{2}$, on construit ainsi des hamiltoniens arbitrairement proche de $H$, sans aucun tore invariant. Cette construction semble manifestement en contradiction avec la persistance des tores invariants, que nous avons vue au n° précédent. 

Pour éclaircir ce paradoxe apparent, considérons la famille dépendant du paramètre $t\in~[0,1]$ définie par
$$H_t=H+t(\frac{p}{q}-\sqrt{2})I_2+\cos(p\p_1-q\p_2). $$
La variante à paramètre du théorème sus-cité montre que  $H_t$ possède de nombreux tores KAM pour presque tout $t$ suffisamment petit, mais ce théorème n'est plus valable pour la valeur $t=1$. Ce qui montre qu'il n'y pas de contradiction.

On trouve, donc dans l'étude des systèmes hamiltoniens, un phénomène  qui n'existe pas pour les singularités d'applications où la stabilité homotopique est équivalente à la stabilité~\cite{Mather_fdet}. 

 On voit, sur cet exemple, que les conditions de non-dégénérescence s'interprêtent comme des conditions d'uniformité sur la taille des voisinages, en fonctions de perturbations.  Pour de nombreux problèmes de physique mathématique, la perturbation est fixée par le contexte, il n'y a alors pas lieu de rechercher une telle uniformité.

\parag Passons à présent, à l’étude des dégénérescences des tores KAM.  Pour cela , considérons à nouveau la variété symplectique $\RM^{2n}$, munie de coordonnées $q = (q_1,...,q_n),\ p=(p_1,...,p_n)$ et de la forme symplectique $\sum_{i=1}^n dq_i \w dp_i$. Rappelons qu'un point critique d'un hamiltonien est
elliptique si le flot de sa partie quadratique est un sous-groupe à un paramètre de rotations. Il est dit {\em diophantien} si la vitesse angulaire de cette famille de rotations, i.e. l'élément de l'algèbre de Lie qui lui correspond, définit un vecteur diophantien.
\begin{theoreme} Le voisinage d'un point critique elliptique diophantien possède un ensemble de mesure positive de tores KAM. De plus, la densité de ces tores tend vers un lorsque l'on approche le point critique.
\end{theoreme} 
Ce théorème, qui répond également à une question de Herman, est en fait un résultat de type KAM. En effet, on peut démontrer que le point critique est l'unique point réel d'une variété lagrangienne complexe singulière, invariante par le flot de $H$. Les tores KAM forment donc une famille qui
dégénère en une variété lagrangienne singulière. Inversement, l'existence d'une telle variété lagrangienne entraîne, comme dans le théorème KAM classique, l'existence d'un ensemble de mesure positive de telles variétés.

Ce deuxième théorème admet une variante  au voisinage d'une orbite périodique. Ainsi, en interprétant tout symplectomorphisme, comme l'application de premier retour de Poincaré d'un champ de vecteur hamiltonien non-autonome, on résout la conjecture de Herman à savoir que: {\em tout symplectomorphisme analytique possède au voisinage d'un point fixe elliptique diophantien, un ensemble de mesure positive de tores invariants}~\cite{Herman_ICM}.

\parag Pour terminer, signalons que la condition diophantienne, peut-être considérablement affaiblie. Pour chaque vecteur $\a \in \RM^n$, définissons la suite  $$\s(\a)_n:=\inf \{ \left| (\a,i) \right|:i \in \ZM^n, \| i \| \leq 2^n \} $$
 où $(\cdot,\cdot)$ désigne le produit scalaire et $\left|\, \cdot \, \right|$ la valeur absolue. Les théorèmes de tores invariants, annoncés dans cette note, restent valable pourvu que la suite $\s(a)$ vérifie la condition arithmétique de Bruno~\cite{Brjuno}. De ce point de vue, un vecteur diophantien définit une suite géométrique, dont la décroissance est très lente comparée à des suites à décroissance exponentielle. Dans le cas général, les tores invariants forment  alors une famille $C^\infty$ tandis que, dans le cas diophantien, elles forment même une famille Gevrey, comme cela avait déjà été montré par Popov pour le théorème classique~\cite{Popov_KAM}.
 
 \noindent{\em Remerciements.} { Ce travail a été financé, en partie, par le  Deutsche Forschungsgemeinschaft project, SFB-TR 45, M086, {\em Lagrangian geometry of integrable systems} et par le Max Planck Institut für Mathematik de Bonn.}
   \bibliographystyle{amsplain}
\bibliography{master}
\end{document}